\documentclass{article}
\usepackage{graphicx} 
\usepackage{xcolor}

\usepackage{amsmath,amsfonts,amsthm,amscd,amssymb,graphicx}

\usepackage{amssymb}

\def\rit{{\Bbb R}}

\def\eps{\varepsilon}

\def\og{\leavevmode\raise.3ex\hbox{$\scriptscriptstyle\langle\!\langle$~}}
\def\fg{\leavevmode\raise.3ex\hbox{~$\!\scriptscriptstyle\,\rangle\!\rangle$}}
\def\beq{\begin{equation}}
\def\eeq{\end{equation}}

\newtheorem{theorem}{Theorem}[section]
\newtheorem{lemma}[theorem]{Lemma}
\newtheorem{e-proposition}[theorem]{Proposition}

\newtheorem{e-definition}[theorem]{Definition\rm}
\newtheorem{remark}{\it Remark\/}

\newtheorem{theoreme}{Th\'eor\`eme}[section]

\newtheorem{proposition}[theoreme]{Proposition}

\title{  Landau damping in mixed hyperbolic-kinetic systems and thick sprays}
\author{D. Bian, B. Despr\'es, V. Fournet, E. Grenier}

\begin{document}

\maketitle


\subsubsection*{Abstract}


This article is devoted to the study of a model of thick sprays which combines the Vlasov equation for the particles and the
barotropic compressible Euler equations to describe the fluid, coupled through the gradient of the pressure of the fluid.
We prove that sound waves interact with particles of nearby velocities,
which results in a damping or an amplification of these sound waves, depending on the sign of the derivative of the distribution
function at the sound speed. This mechanism is very similar to the classical Landau damping which occurs in the Vlasov-Poisson
system. 
If the sound waves are amplified then the thick spray model is linearly ill-posed in Sobolev spaces, even locally in time.

We also show that such Landau damping type phenomena naturally arise
when we couple an hyperbolic system of conservation laws with the Vlasov equation.


\section{Introduction}


We consider a compressible fluid described by its volume fraction $\alpha(t,x)$,
its density $\rho(t,x)$ and its velocity $u(t,x)$. We moreover assume that this fluid
is ideal and barotropic, with a pressure law $p(\rho)$.
This fluid contains particles which are small spheres, of radius $r_p$,
and which are described by a distribution function $f(t,x,v)$.
We will assume that the only effect of the particles is that they occupy some fraction $\alpha(t,x)$ of the available volume,
 given by
\beq \label{thick4}
\alpha(t,x) = 1 - \kappa \int_{\rit^3} f(t,x,v) \, dv
\eeq
where 
 the parameter $\kappa$ is defined by
\begin{equation} \label{eq:bd3}
\kappa = {4 \over 3} {\pi r_p^3 \over V_{\rm ref}}>0,
\end{equation}
where $r_p$ is the radius of a particle and $V_{\rm ref}$
is a reference volume. 
This leads to the following non dimensional barotropic compressible Euler equations for the fluid part
\beq \label{thick1}
\partial_t (\alpha \rho) + \nabla_x \cdot (\alpha \rho u) = 0,
\eeq
\beq \label{thick2}
\partial_t (\alpha \rho u) + \nabla_x \cdot ( \alpha \rho u \otimes u) + \alpha \nabla_x p(\rho) = 0.
\eeq
We will also assume that the particles react to the pressure gradient $\nabla_x p(\rho)$
of the fluid, which leads to
\beq \label{thick3}
\partial_t f + v \cdot \nabla_x f - \nabla_x p(\rho) \cdot \nabla_v f = 0.
\eeq
In particular, we neglect the collisions between the particles.
The equations (\ref{thick4},\ref{thick1},\ref{thick2},\ref{thick3}) are a set
of four equations on the four unknowns $\alpha$, $\rho$, $u$ and $f$, 
the first one being a  diagnostic variable, 
the last  three of them being prognostic variables.

This system is a prototype for describing so-called ``thick" sprays. Such coupling is used when the particles are small but occupy a non-negligible volume fraction of the mixture~\cite{williams,orourke,dukowicz}.
A linear friction force of the form $\beta (u - v)$ where $v$ is the  velocity of the particles and $\beta \in \rit$
is the friction parameter is usually present,
but we disregard this term for the simplicity of
the mathematical exposure.
We will restrict ourselves to one space dimension, the computations being similar in three space dimensions.

\medskip

Up to the best of our knowledge, the existence of solutions  to this system  with Sobolev regularity, even locally in time is an open question. In this note, we study the linearized system in the asymptotic regime, namely the ``thin spray" regime ($\kappa \to 0$). In this regime, we study in detail the Landau damping type phenomenon attached to this system, discovered by C. Buet, B. Despr\'es and V. Fournet in \cite{Buet,Fournet}. Two cases arise: either the acoustic waves are amplified or they are damped. In the first case,
we will prove that (\ref{thick4},\ref{thick1},\ref{thick2},\ref{thick3}) is linearly  ill-posed  in Sobolev spaces thereby
confirming a conjecture of C. Baranger and L. Desvillettes \cite{Baranger}. 
In the second case, we 
conjecture that  (\ref{thick4},\ref{thick1},\ref{thick2},\ref{thick3}) is well-posed in Sobolev spaces, locally in time.

\medskip

Let us first describe the linearized equations around a constant state.
Let $0 < \alpha_0 < 1$, let $\rho_0 > 0$, and let $f_0(v)$ be a given distribution function.
Up to a Galilean change of variables, and the corresponding change of $f_0$, we may assume that $u_0 = 0$.
We naturally  assume the compatibility condition
\beq \label{eq:compi}
\alpha_0 = 1 - \kappa\int_{\rit^3} f_0(v) \, dv .
\end{equation}
We note that (\ref{thick1},\ref{thick2}) may be combined in
$$
\rho \partial_t u + \rho u \cdot \nabla_x u + \nabla_x p(\rho) = 0.
$$
The linearized system around ($\alpha_0$, $\rho_0$, $u_0 = 0$, $f_0$) is thus 
\beq \label{thick1l}
\alpha_0 \partial_t \rho_1 + \alpha_0 \rho_0 \nabla_x \cdot u_1 
+  \kappa \rho_0 \nabla_x \cdot \int_{\rit^3} f_1 v \, dv = 0,
\eeq
\beq \label{thick2l}
\rho_0 \partial_t u_1 +  p'(\rho_0)\nabla_x \rho_1 = 0,
\eeq
\beq \label{thick3l}
\partial_t f_1 + v \cdot \nabla_x f_1 - p'(\rho_0) \nabla_x \rho_1 \cdot \nabla_v f_0  = 0.
\eeq
We introduce the speed of sound in the fluid (without particles)
$$
c_0 = \sqrt{p'(\rho_0)}.
$$
We note that $\kappa$ goes to $0$ when $r_p$ goes to $0$, namely when the volume occupied
by the particles goes to $0$. The limit $\kappa \to 0$ thus corresponds to the limit of ``thin sprays". 

\medskip

We will say that the linear system (\ref{thick1l}-\ref{thick3l}) is spectrally stable (respectively spectrally unstable) if it has no solution of the form
\begin{equation} \label{eq:vec}
(\alpha_1,\rho_1,u_1,f_1) = (\alpha_\star,\rho_\star,u_\star,f_\star) e^{i k \cdot x - i \omega t}
\end{equation}
with $\Im \omega > 0$ (respectively if it has one solution of this form with $\Im \omega > 0$).

\begin{theorem} \label{theo0}
Let $\alpha_0 \in (0,1)$, $\rho_0 > 0$, $u_0 \in \rit$ and $f_0(v)$ be given. 
Let us assume that $f_0(v)$ is analytic and can be extended to the complex strip $| \Im v | \le \delta$ for some positive $\delta$.
Let us moreover assume that, on this strip, $| f(v) | \le C_0 \exp(- C_1 |v|^2)$ for some constants $C_0$ and $C_1$.
Then,
\begin{itemize}
\item if $f_0(v) = \mu(|v - u_0|^2)$ where $\mu$ is a smooth and decreasing function,  (\ref{thick1l}-\ref{thick3l}) is spectrally stable,
    
\item if $\kappa$ is small enough, (\ref{thick1l}-\ref{thick3l}) is spectrally stable if  $\partial_v f(u_0 \pm c_0^\star) < 0$
and  spectrally unstable if  $\partial_v f(u_0 \pm c_0^\star) > 0$, where $c_0^\star$ is the sound speed in the spray.
   
\item If (\ref{thick1l}-\ref{thick3l}) is spectrally unstable, then (\ref{thick4},\ref{thick1},\ref{thick2},\ref{thick3}) is linearly
ill-posed in Sobolev spaces.

\end{itemize}

\end{theorem}

The first point of this theorem has already been proved in \cite{JSP}.
The second point is proved in section \ref{sec:kappa} and the third point in section \ref{sec:unstable}.

In the last part of this work, we prove Proposition
\ref{pro:1} which is a necessary condition for observing a  linear  Landau damping arising in  kinetic perturbations of quasi-linear strictly  hyperbolic systems 
of conservation laws.

The value of $c_0^\star$, sound speed in the spray (fluid
with particle), is provided in equation (\ref{velo}).
When $\kappa = 0$, equations
(\ref{thick1l},\ref{thick2l}) are  decoupled from the kinetic part
and are the linearized barotropic compressible Euler equations.
They admit traveling wave solutions, namely the classical acoustic waves, of speed $c_0$.
When $\kappa \ll 1$, we will prove that the system (\ref{thick1l}-\ref{thick3l}) has waves with a phase velocity
$c_0^\star = c_0 + O(\kappa)$. These waves are  damped provided $\partial_v f(u_0 \pm c_0^\star)< 0$ 
but they are amplified if $\partial_v f(u_0 \pm c_0^\star) > 0$.


Our proofs are actually more general than just studying possible unstable modes like 
(\ref{eq:vec}).
Indeed we  construct the dispersion relation 
 $$
 \omega \mapsto \mathcal D(\omega,k)\in \mathbb C
 \mbox{ for }\omega \in \mathbb C
 $$
 associated 
to the linearized equations 
(\ref{thick1l}-\ref{thick3l}).
A  comprehensive  mathematical study of the 
stability property of the linearized 
Vlasov-Poisson equations based on the dispersion relation is 
Degond's seminal work \cite{Degond}.
Some conclusions can be generalized to our case. In particular it explains the structure of the solutions associated 
to roots $\omega$ of the dispersion relation with negative imaginary parts $\Im \omega<0$.
Some  connections between the existence of unstable modes and the stability of the Vlasov-Benney system are in \cite{Bardos}.

\medskip 

If we add a viscosity to (\ref{thick3l}), namely if we describe the fluid part by the compressible
Navier-Stokes equations, all the waves with a sufficiently large wave number $k$ are damped. As proved in \cite{HanKwan}, if moreover
we assume a Penrose type assumption on $f_0(v)$, the corresponding system is locally well-posed in Sobolev spaces.

\medskip

This note is organized as follows: 
the second section is devoted to the study of the dispersion relation of the linearized system
(\ref{thick1l}-\ref{thick3l}), to the study of the limit $\kappa \to 0$, and to the comparison with the classical Landau  damping. 
The third section is devoted to the proof of the ill-posedness in the spectrally unstable case.
In the last section we show that Landau damping naturally occurs when we
couple an hyperbolic system of conservation laws and the Vlasov equation and we provide a characterization of stable systems in Proposition \ref{pro:1}.


\section{Study of the dispersion relation}



\subsection{The dispersion relation}


The aim of this section is to compute the dispersion relation of (\ref{thick4},\ref{thick1},\ref{thick2},\ref{thick3}).
We assume that  $\alpha_0$, $\rho_0$, $u_0$ and $f_0(v)$ are given.
Up to a Galilean change of variables, we may linearize (\ref{thick4},\ref{thick1},\ref{thick2},\ref{thick3}) around $u_0 = 0$, up to a translation of $f$. After this change of variables,
the velocity $v$ is shifted by $-u_0$ and the initial density $f_0$ is redefined accordingly. 
We also  assume that $f_0(v)$ is analytic and can be extended to the complex strip $| \Im v | \le \delta$ for some positive $\delta$.

\begin{lemma}
Let $\alpha_0$, $\rho_0$, $u_0 = 0$ and $f_0(v)$ be given.
Let $k \in \mathbb{Z}^*$ or $\mathbb{R}^*$. 
    The  dispersion relation ${\cal D}(k,\omega)$ for $\omega\in \mathbb C^*$ only depends on 
    $$
    \sigma = {\omega \over |k|} 
    $$
    and is, for $\Im \sigma > 0$,
    \begin{equation}  \label{dispersion}
    {\cal D} =  1 -  {c_0^2 \over \sigma^2} 
- \frac{\kappa \rho_0 c_0^2}{\alpha_0} {1 \over \sigma} \int_{\rit^3} { v \partial_{v} f_0 \over v - \sigma}  \, dv .
\end{equation}
It can be extended in an holomorphic way to $\Im \sigma > - \delta$
through the formula
 \begin{equation}  \label{dispersion}
    {\cal D} =  1 -  {c_0^2 \over \sigma^2} 
- \frac{\kappa \rho_0 c_0^2}{\alpha_0} {1 \over \sigma} P.V. \int_{\rit^3} { v \partial_{v} f_0 \over v - \sigma}  \, dv 
- i \pi \kappa \rho_0 c_0^2 \Bigl[ \mathbf 1_{\Im \sigma = 0} + 2  \mathbf 1_{\Im \sigma < 0} \Bigr] \partial_v f_0 ( \sigma ) 
\end{equation}
where $P.V.$ denotes the principal value and  $\mathbf 1_{\dots}$ denotes
the indicatrix function.
\\
Moreover 
if $f_0$ is a monotonic function of $|v|^2$, the dispersion relation has no root $\sigma$ with $\Im \sigma > 0$.
\end{lemma}

\begin{remark}
The dispersion relation (\ref{dispersion}) makes sense for
$\Im \sigma \ge -\delta$. 
 Also the term between brackets can be rewritten as
$\mathbf 1_{\Im \sigma = 0} + 2  \cdot \mathbf 1_{\Im \sigma < 0}=
1-\mathrm{sign} \, \Im( \sigma )$. 
\end{remark}

\begin{remark}
Note that if $\sigma$, with $\Im \sigma < 0$, is a zero of the extended dispersion relation, it is not a ``regular" eigenvalue of the thick spray model, but only a decay rate of the linearized system (see \cite{Degond} for a discussion in the case of the Vlasov-Poisson system). 
\end{remark}

\begin{proof} 
Let
$
\tau_1 = - \rho_0^{-2} \rho_1$ be the linearization of the specific volume.
Then
\beq \label{thick1ll}
\alpha_0 \rho_0 \partial_t \tau_1 = \alpha_0 \nabla_x \cdot u_1 
+  \kappa  \nabla_x \cdot \int_{\rit^3} f_1 v  \, dv,
\eeq
\beq \label{thick2ll}
 \rho_0 \partial_t u_1 = \rho_0^2 c_0^2 \nabla_x \tau_1,
\eeq
\beq \label{thick3ll}
\partial_t f_1 + v \cdot \nabla_x f_1 + c_0^2 \rho_0^2 \nabla_x \tau_1 \cdot \nabla_v f_0  = 0.
\eeq
We take the Fourier-Laplace transform of this system, with dual Fourier variables $k$ in space and $- \omega$ in time where $\Im\omega>0$, which leads to
\beq \label{tt}
-i \omega   \rho_0 \widehat \tau_1 =    i k \cdot \widehat u_1 +   \frac \kappa{\alpha_0} i k \cdot \int_{\rit^3} \widehat f_1 v  \, dv +  \rho_0\widehat \tau^{init} 
\ ,
\eeq
\beq
-i \omega  \rho_0 \widehat u_1 = i k \rho_0^2 c_0^2 \widehat \tau_1 +  \rho_0 \widehat u^{init},
\eeq
\beq
( - i \omega  + i k \cdot v )  \widehat f_1 =  - \rho_0^2 c_0^2 i k \cdot \nabla_v f_0 \widehat \tau_1 + \widehat f^{init}.
\eeq
Thus
$$
\widehat u_1 = - {k \over \omega} {\rho_0 c_0^2 } \widehat \tau_1 + {i \over \omega} \widehat u^{init}
$$
and
$$
\widehat f_1 = - \rho_0^2 c_0^2 { k \cdot \nabla_v f_0 \over k \cdot v - \omega} \widehat \tau_1
+ {\widehat f^{init} \over i k \cdot v - i \omega},
$$
which together with (\ref{tt}),  gives
\begin{align*} 
-i \omega   \rho_0 \widehat \tau_1 &= i  k \cdot
\left( - {k \over \omega} \rho_0 c_0^2  \widehat \tau_1 + {i \over \omega} \widehat u^{init}\right)
\\
& \quad +\frac\kappa{\alpha_0} i k \cdot \int_{\rit^3} 
\left(- \rho_0^2 c_0^2 { k \cdot \nabla_v f_0 \over k \cdot v - \omega} \widehat \tau_1
+ {\widehat f^{init} \over i k \cdot v - i \omega} \right)
v  \, dv
 +  \rho_0\widehat \tau^{init} .
\end{align*}
This is reorganized  as 
$$
\begin{array}{lll}
0 & =& \displaystyle \left(1-c_0^2 \frac{k^2}{\omega^2} 
-
 {\frac\kappa{\alpha_0}  \rho_0 c_0^2} {k \over \omega} \int_{\rit} {\partial_{v} f_0 \over v - {\omega \over k}} v\, dv\right)i \omega    \widehat \tau_1 \\
 & & \displaystyle+\tau^{init} -\frac{k\cdot \widehat u^{init}}{\omega \rho_0}
 +\frac{\kappa}{\rho_0\alpha_0} \int_{\rit} {\widehat f^{init}(v) \over v - {\omega \over k}} v \, dv.
 \end{array}
$$
This expression is rewritten as
\beq \label{dispersion0}
{\cal D}_+(k,\omega) \Bigl[ -  i {\omega } \widehat \tau_1 \Bigr] = {\cal E}(k,\omega) 
\eeq
where
\beq \label{dispersion1}
{\cal D}_+(k,\omega) = 1 - c_0^2 {k^2 \over \omega^2} 
- \frac{\kappa \rho_0 c_0^2}{\alpha_0} {k \over \omega} \int_{\rit} {\partial_{v} f_0 \over v - {\omega \over k}} v\, dv
\eeq
and where the right-hand side is 
\beq \label{dispersion2}
{\cal E}(k,\omega) = {1 \over \rho_0 } \Bigl[
  \rho_0 \widehat \tau^{init}  - {k \cdot \widehat u^{init} \over \omega} + \frac\kappa{\alpha_0} \int_{\rit} {\widehat f^{init}(v) \over v 
  -{\omega \over k}} v \, dv
\Bigr] .
\eeq
The expression (\ref{dispersion0}) has been derived  for
$\Im \omega > 0$. Since the   integral in 
(\ref{dispersion}) is singular for $\omega = k v$, then (\ref{dispersion0}) is clearly meaningless for 
$\Im \omega \leq  0$.

 Following Landau's approach, we extend ${\cal D}_+$  by analytic continuation in the complex plane.
 We define, for $\Im \omega = 0$,
$$
{\cal D}_0(k,\omega) = 1 - c_0^2  {k^2 \over \omega^2} 
- \frac{\kappa \rho_0 c_0^2}{\alpha_0} {k \over \omega}  P.V. \int_{\rit} { \partial_{v} f_0 \over v - {\omega \over k}} v \, dv - i \pi \kappa \rho_0 c_0^2 \partial_v f_0 \Bigl( {\omega \over k} \Bigr)
$$ 
where $P.V.$ denotes a principal value and, for $\Im \omega < 0$,
$$
{\cal D}_- (k,\omega) = 1 - c_0^2 {k^2 \over \omega^2} 
-  \frac{\kappa \rho_0 c_0^2}{\alpha_0} {k \over \omega}  \int_{\rit} { \partial_{v} f_0 \over v - {\omega \over k}} v \, dv 
- 2 i \pi \kappa \rho_0 c_0^2 \partial_v f_0 \Bigl( {\omega \over k} \Bigr).
$$
We will denote by ${\cal D}(k,\omega)$ the full relation dispersion, which equals
${\cal D}_+$ if $\Im \omega > 0$, ${\cal D}_0$ if $\Im \omega = 0$ and
${\cal D}_-$ if $\Im \omega < 0$.
We note that ${\cal D}(k,\omega)$ only depends on $\sigma = \omega / |k|$.

Now let us assume that $\sigma$ is an unstable eigenvalue, in the sense $\Im \sigma > 0$,  then
$$
\sigma - {c_0^2 \over \sigma} 
- \frac{\kappa \rho_0 c_0^2}{\alpha_0}  \int_{\rit} { v \partial_{v} f_0 \over v - \sigma}  \, dv = 0.
$$
Thus
$$
\sigma - c_0^2  {\bar \sigma \over |\sigma|^2} 
- \frac{\kappa \rho_0 c_0^2}{\alpha_0}  \int_{\rit} { v \partial_v f_0 \over |v - \sigma|^2} (v - \bar \sigma)  \, dv = 0.
$$
The imaginary part gives
$$
\Bigl[ 1 +  {c_0^2 \over |\sigma|^2} 
- \frac{\kappa \rho_0 c_0^2}{\alpha_0}  \int_{\rit} { v\partial_{v} f_0 \over |v - \sigma|^2}   \, dv \Bigr] \ \Im \sigma  = 0.
$$
As $v \partial_v f_0 \le 0$, the quantity between brackets is positive, thus $\Im \sigma = 0$. It is a  contradiction so it  ends the proof of the last part of the lemma.
Note that this proof is very similar to the proof of the classical Rayleigh's criterium in fluid mechanics \cite{Drazin}.
\end{proof}


\subsection{Comparison with the genuine Landau damping}


Let us now compare the previous dispersion law with the dispersion law of the genuine Landau damping \cite{Landau,Schekochihin}. 
The classical Vlasov-Poisson system takes the form of 
$$
\partial_t f + v \cdot \nabla_x f + E \cdot \nabla_v f = 0,
$$
$$
E = - \nabla_x V, \qquad
- \Delta V = \int_{\rit^3} f \, dv,
$$
where $f(t,x,v)$ is the distribution function, $E(t,x)$ the electrostatic field and $V(t,x)$ the electrostatic potential. In this case, the dispersion relation is in one space dimension
$$
{\cal D}_{Landau}(k,\omega) = 1 
- {1 \over |k|^2} \int_{\rit} {\partial_{v} f_0 \over v - {\omega \over k}} \, dv + i \pi 
\Bigl[ \mathbf 1_{\Im \omega = 0} + 2 \cdot \mathbf 1_{\Im \omega < 0}\Bigr] \partial_v f \Bigl( {\omega \over k} \Bigr).
$$
Thus, the two dispersion relations ${\cal D}$ and ${\cal D}_{Landau}$ have a very close structure, except that the dispersion relation for thick spray only depends on $\omega / |k|$, whereas the Landau dispersion relation depends on $\omega / k$ and $k$.

This difference comes from the fact that the force in the Vlasov-Poisson system is of order $-1$ whereas in our case, $\nabla_x p(\rho)$ is of order $+1$, like the transport part $\partial_t + v \cdot \nabla_x$. In the Vlasov-Poisson case, when $|k| \gg 1$, the transport term is dominant and the electric field plays a negligible role, whereas, when $|k| \ll 1$, it is dominant. In the physical space, the behavior of the Vlasov-Poisson system is
completely different on scales much smaller than the Debye length (where the transport is predominant) and on scales much larger
than this length (where the electric field is predominant). 

In our system on the contrary, all the terms of (\ref{thick1},\ref{thick2},\ref{thick3}) are of order $+1$ and scale in the same way. As a consequence, the
linearized system only depends on the phase speed parameter $\omega / k$.

This remark is particularly important if one unstable mode $\sigma$ is found, with $\Im \sigma > 0$, since the corresponding waves
satisfy $\omega = \sigma k$, which means that the time scale of instability goes to $0$ as $|k|$ goes to infinity.


\subsection{The ``thin spray" regime \label{smallkappa} \label{sec:kappa}}


The dispersion relation (\ref{dispersion}) depends on two physical parameters which are $c_0$, the speed of sound in the fluid, 
and $\kappa$, which is correlated to the radius of the particles.

To prepare the next expansions, we write $\mathcal D = \mathcal D_r + i \mathcal D_i$, where
$$
\mathcal D_r(\sigma)=   1 -  {c_0^2 \over \sigma^2} 
- \frac{\kappa \rho_0 c_0^2}{\alpha_0} {1 \over \sigma} P.V. \int_{\rit} { v \partial_{v} f_0 \over v - \sigma}  \, dv 
$$
and 
$$
\mathcal D_i(\sigma)=
 -\pi \kappa \rho_0 c_0^2 \Bigl[ \mathbf 1_{\Im \sigma = 0} + 2  \mathbf 1_{\Im \sigma < 0} \Bigr] \partial_v f_0 ( \sigma ).
$$
Note that $\mathcal D_r$ is defined for all $\sigma \in \mathbb C^*$.
The other term $\mathcal D_i$ is restricted to the strip 
$|\Im \sigma| \leq \delta$.

Expansions at various orders $O(\sigma^{-n})$ of the first part are easily obtained as follows. Let
\begin{equation} \label{eq:bd222}
{\cal F}(\sigma) = {1 \over \sigma} P.V. \int_{\rit} {\partial_{v} f_0 \over v - \sigma} v\, dv .
\end{equation}
For large $|\sigma|$, we expand ${\cal F}(\sigma)$ in $\sigma$ 
\begin{align}
{\cal F}(\sigma) &= - {1 \over \sigma^2} P.V. \int_\rit {\partial_{v} f_0 \over  1 - {v \over \sigma}} v \, dv
\\
&=  - {1 \over \sigma} \int_\rit \sum_{j = 0}^N {v^{j+1} \over \sigma^{j+1}} \partial_ {v} f_0  \, dv
- {1 \over \sigma} P.V. \int_\rit {v^{N+2} \over \sigma^{N+2}} { \partial_{v} f_0 \over 1 - {v \over \sigma}} \, dv.
\end{align}
The integrals with odd $j$ vanish by symmetry.
Integrals with even $j$ are evaluated by integration by parts.
It  gives
for example $$
{\cal F}(\sigma) = {m_0 \over \sigma^2} + {3 m_2 \over \sigma^4} + O(\sigma^{-6})
$$
where
\begin{equation} \label{eq:bd2}
m_0 = \int_\rit f_0(v) \, dv, \qquad
m_2 = \int_\rit f_0(v) v^2 \, dv .
\end{equation}
Thus as $|\sigma| \to + \infty$, one can write 
\beq \label{asymptot}
{\cal D}_r(\sigma) = 1 - \frac{c_0^2}{ \sigma^2} - 
\frac{\kappa \rho_0 c_0^2}{\alpha_0}
{m_0 \over \sigma^2}
- \frac{\kappa \rho_0 c_0^2}{\alpha_0} {m_2 \over \sigma^4} + O(\sigma^{-6}).
\eeq
In particular, ${\cal D}_r(\sigma)$ converges to $1$ as $| \sigma |$ goes to $+ \infty$.

The imaginary part $D_i(\sigma)$ is proportional to  $\partial_v f_0(\sigma)$. It is physically natural to assume that $f_0$ and its derivative
are dominated by Maxwellians at infinity in a strip in the complex plane. So  one has in the strip 
\begin{equation}
\left| f_0(v)\right|+
\left| \partial_v f_0(v)\right|  = O(|v|^{-n})
\mbox{ for }|\Im v|\leq   \delta  \mbox{ (for all } n\geq 0).
\end{equation}
One obtains as well in the strip
\begin{equation} \label{eq:bd111}
{\cal D}_i(\sigma) =  O(|\sigma|^{-n})
\mbox{ for }|\sigma|\to \infty.
\end{equation}

\begin{remark}
The combination of  (\ref{asymptot}) and (\ref{eq:bd111}) implies that ${\cal D}$ has no large root in the strip.
\end{remark}






Let us now turn to the study of the limit $\kappa \to 0$ which physically corresponds  to the fact that the radius of the particles tends to zero 
$r_p\to 0$.
We start from (\ref{dispersion}) and we assume that all terms that appear  are  fixed, except $\kappa$, which goes to $0$. 
In particular $m_0$ is fixed, so (\ref{eq:bd3}) and (\ref{eq:bd2}) imply that $\alpha_0\to 1$.
We obtain
$$
\lim_{\kappa \to 0} 
\sigma^2 {\cal D} = \sigma^2  - c_0^2  ,
$$
thus, when $\kappa = 0$, there are only two roots $\sigma_\pm$, given by
$$
\sigma_\pm = \pm c_0,
$$
which is physically expected: when $\kappa = 0$, namely when there
is no particle in the fluid, the only waves which propagate are the sound
waves, with speed $\pm c_0$.

Now, using (\ref{asymptot}), we see that, provided $\kappa$ is small
enough, there is no zero in the area $|\sigma | \ge 2 c_0$.
As a consequence, using the implicit function theorem, if $\kappa$ is small enough, the dispersion relation has only two roots $\sigma_\pm(\kappa)$, close to $\sigma_\pm$.
To precise the behavior of $\sigma_\pm$ as $\kappa$ goes to $0$, we write
$$
{\cal D}(\sigma) = 0 \approx {\cal D}(c_0) + {\cal D}'(c_0) ( \sigma - c_0)
$$
which gives 
$$
\sigma = c_0 \Bigl[ 1 + {\kappa \over 2} {P.V} \int_\rit {\partial_v f_0 \over v - c_0} \, dv 
+  i \pi \kappa \partial_v f_0(c_0) \Bigr] + O(\kappa^2). 
$$
In particular, the sound speed in the spray is
\begin{equation} \label{velo}
c_0^\star =  c_0 \Bigl[ 1 + {\kappa \over 2} {P.V} \int_\rit {\partial_v f_0 \over v - c_0} \, dv \Bigr]  +   O(\kappa^2). 
\end{equation}
To obtain an equivalent of the imaginary part of $\sigma$ we write
(without more justification) that the dispersion identity 
${\cal D}(\sigma) = 0$ can be approximated in the strip 
 by
$$
 {\cal D}_r(c_0^\star ) + {\cal D}_r'(c_0^\star ) (\sigma -c_0^\star )
 + i \mathcal D_i(c_0^\star )=0.
 $$ 
Since ${\cal D}_r(\sigma_0)=0$ by construction, 
one obtains 
\beq \label{damping}
\Im \sigma \approx - i  { {\cal D}_i(c_0^\star ) \over {\cal D}'_r(c_0^\star )},
\eeq
where
$$
\mathcal D_i(c_0^\star )= - \pi \kappa \rho_0 c_0^2  \partial_v f_0 ( c_0^\star  ).
$$
Note that the corresponding wave is damped or amplified, depending on the sign of $\partial_v f(c_0)$.

Let us now discuss (\ref{velo}). 
Particles with a speed close to $c_0^\star$ see almost no change in these waves and undergo large variations in their velocities. 

In the genuine Landau damping, particles with a velocity close to the phase velocity of the electric waves strongly interact with the electric waves. 
Particles with slightly smaller velocities are accelerated by the electric wave and take energy from the wave. 
On the contrary, particles with slightly larger velocities are decelerated and give part of their kinetic energy to the wave.

As a consequence, if $f_0$ is decaying near the phase velocity of the electric field,
more particles are accelerated than decelerated. The net effect is a transfer of the energy
of the electric field to the kinetic energy of the particles: the electric wave is damped.

The situation is similar here. Waves propagate in the thick spray with a velocity $c_0^\star$. 
 Slightly slower particles  get accelerated and slightly faster particles  get decelerated.
If $\partial_v f_0(c_0^\star) < 0$, there is a net transfer of energy from the waves to the particles.
As a consequence, the waves are damped.

If on the contrary $\partial_v f_0(c_0^\star) > 0$, there is a transfer of energy from the particles
to the waves, which are thus amplified, leading to an instability.


\section{The spectrally unstable case} \label{sec:unstable}


We now prove that the thick spray system is ill-posed if the
Landau damping is negative, namely if there exists an eigenvalue
with a positive imaginary part (last statement of Theorem \ref{theo0}).

The proof is classical and is just a play with scalings.
Let $\omega$ be a root of the dispersion relation with
with $\Im \omega > 0$.
Then the  linearized system (\ref{thick1l}-\ref{thick2l}) or
equivalently the linearized system
(\ref{thick1ll}-\ref{thick3ll}) has  solutions 
whose $L^2$ norm grow like $e^{\Im \omega t}$.

The verification for 
(\ref{thick1ll}-\ref{thick3ll}) is easy.
Assume that $\omega$ with $\Im \omega >0$ is a  root of the dispersion
relation $\mathcal D(\omega,k)=0$ and let
\begin{equation} \label{eq:fofo}
\left\{
\begin{array}{rll}
&\tau(t,x) =  \Re\left(e^{i k x - i \omega k t}  \right), \\
&u(t,x) =  - k \rho_0 c_0^2\  \Re\left(e^{i k x - i \omega k t}  \right), \\
&f(t,x,v) = -\rho_0^2 c_0^2 f_0'(v) \
 \Re\left( \frac{\displaystyle e^{i k x - i \omega k t} }{v-\omega/k} \right).
\end{array}
\right. 
\end{equation}
It is immediate to check the functions (\ref{eq:fofo}) are solutions to the linearized
equations. Then it is generalized to (\ref{thick1l}-\ref{thick2l})
after an evident change of variables.

Let us now concentrate on the consequences for (\ref{thick1l}-\ref{thick2l}).
Let $s \ge 0$ be fixed and consider the linearized system (\ref{thick1l}-\ref{thick2l}). Let $N > s$ and let
$$
(\rho_k,u_k) = k^{-N} \Re \left( 
(\rho_\star,u_\star) e^{i k x - i \omega k t} 
\right).
$$
Then, at $t = 0$, $\| (\rho_k,u_k) \|_{H^s} \to 0$, and at $t_k = (N+1) k^{-1} \log k$,
$\| (\rho_k,u_k) \|_{L^2} \ge \theta_0 > 0$ for some positive $\theta_0$ as $k \to + \infty$.
Note that $t_k \to 0$ as $k \to + \infty$.
The linearized system is thus ill-posed in Sobolev spaces. This argument ends the proof of Theorem 
\ref{theo0}.

We now give a bump-on-tail  example of 
an unstable profile.
More precisely we prove that if $f_0$ is a given smooth distribution function,
decaying in $|v|$, then we can find some arbitrarily small perturbation
$g$ such that $f = f_0 + g$ is spectrally unstable.
Let $g(v)$ be a smooth, non negative function, with support in $[-1,+1]$, such that $\partial_v g(0) > 0$ and such that
$$
\int_\rit g(v) \, dv = 1.
$$
Let $\varepsilon > 0$ be arbitrarily small. Let $\eta > 0$.
Let $c_0^\star$ be large enough.
We set
$$
f(v) = (1 - \eps) f_0(v) + \eps \eta g \Bigl( {v - c_0^\star \over \eta} \Bigr) \int_\rit f_0(v) \, dv.
$$
We note that $f$ and $f_0$ have the same integral. Moreover, provided $\eta$ is small enough and $c_0^\star$ is large enough, $\partial_v f(c_0^\star) > 0$.
Thus,  the dispersion relation has a root $\sigma$ with $\Im \sigma > 0$
and hence $f$ is spectrally unstable.

A physical interpretation is as follows: particles of the bump slow down, thereby releasing kinetic energy which, by conservation of total energy,  is transferred to the wave which increases.


\section{  ``Universality" of the Landau damping}


Our objective is to show that the dispersion relation is easy 
to obtain for generic linearized  conservation laws
coupled with a kinetic equation.
It follows that   linear Landau damping is universal in the sense that 
it can arises for many physical problems.

\subsection{The scalar case}

Let us first study the coupling of a scalar conservation law with the 
Vlasov equation. We consider an ``hyperbolic coupling", namely
a coupling through derivatives, like for instance
\beq \label{coupl1}
\partial_t u + \lambda(u) \partial_x u + \kappa \partial_x \int f(v) v \, dv = 0,
\eeq
\beq \label{coupl2}
\partial_t f + v \partial_x f + \partial_x u \, \partial_v f = 0,
\eeq
where $\lambda(u)$ is a given function and $\kappa \in \rit$ is a coupling constant, which will be assumed to be small $|\kappa| \ll 1$.

Let $u_0 \in \rit$ and $f_0(v)$ be a given distribution function. Then, 
the linearization of (\ref{coupl1},\ref{coupl2}) gives
\beq \label{coupl1l}
\partial_t \widetilde u + \lambda(u_0) \partial_x \widetilde u + \kappa \partial_x \int \widetilde f(v) v \, dv = 0,
\eeq
\beq \label{coupl2l}
\partial_t \widetilde f + v \partial_x \widetilde f + \partial_x \widetilde u \, \partial_v f_0 = 0.
\eeq
Let us consider an eigenmode $e^{i k (x - \omega t)} (\widetilde u, \widetilde f)$.
Then we obtain
$$
\widetilde f(v) = - \int_\rit {\partial_v f_0(v) \over v - \omega} \, dv
$$
and thus
\beq \label{scalar1}
\omega =  \lambda(u_0) - \kappa \int_\rit {v \partial_v f_0(v)\over v - \omega} \, dv.
\eeq
We must keep in mind that (\ref{scalar1}) is obtained by taking the Laplace transform of (\ref{coupl1},\ref{coupl2}),
namely through a contour integral where $\Im \omega$ is large enough.

This formula may then be extended for $\Im \omega > 0$ through analytic continuation. However we have to take
care that the integrand is singular when $v = \omega$. According to Plemelj's formula, to extend (\ref{scalar1}) 
to real $\omega$, we have to replace the integral by a principal value and to add an imaginary extra term.
For $\Im \omega  = 0$, the dispersion relation reads
\beq \label{scalar2}
\omega =  \lambda(u_0) - \kappa P.V. \int_\rit {v \partial_v f_0(v)  \over v - \omega} \, dv - i \pi \kappa \omega \partial_v f_0(\omega).
\eeq
If $f_0$ is monotonic in $|v|$, namely if $\partial_v f(v) \ne 0$ if $v \ne 0$,
them we see that the only possible real solution of (\ref{scalar2}) is $\omega  = 0$.

Let us assume that $\lambda(u_0) \ne 0$. Then $\omega$ is a smooth function of $\kappa$. By continuity $\omega(\kappa)$ may
not vanish for small $\kappa$. As a consequence, either $\Im \omega(k) < 0$, which corresponds to a damping of the hyperbolic system
through its interaction with the particles, a kind of of Landau damping, or $\Im \omega(k) > 0$, corresponding to an amplification
of the hyperbolic system (negative Landau damping).

Let us study the case $|\kappa| \ll 1$. When $\kappa$ is small, $\omega(\kappa)$ is close to $\lambda(u_0)$, thus
\beq \label{scalar2}
\Im \omega \sim - \pi \kappa \lambda(u_0) \partial_v f \Bigl( \lambda(u_0) \Bigr).
\eeq
If $f_0$ is monotonic in $v$, we get a Landau damping when $\kappa < 0$ and a negative Landau damping when $\kappa > 0$.

\subsection{The case of a system}

Let us now turn to hyperbolic systems of conservation laws, coupled with a Vlasov
equation through derivatives and study systems of the form
\beq \label{sys1}
\left\{
\begin{array}{lll}
\partial_t u + A(u) \partial_x u + \kappa \partial_x \int_\rit f(v) \phi(u,v) \, dv = 0, \\
\partial_t f + v \partial_x f + \partial_x \psi(u) \, \partial_v f = 0,
\end{array}
\right.
\eeq
where $u$ is a vector valued function, with values in $\rit^N$, $A$ is a $N \times N$ matrix, 
$\phi$  is a function from $\rit^N \times \rit$ to $\rit^N$,
$\psi$ is a function from $\rit^N$ to $\rit$ and $\kappa$ is a coupling constant. 
The term $\phi(u,v)$ is responsible of the   "action" of the 
particles on the hyperbolic part.
The term 
$ \partial_x \psi(u)$ represents the "force" that the hyperbolic 
part exerts on the particles.

Let $u_0 \in \rit^N$ be a space-time constant state and
let $v\mapsto f_0(v) \in \rit$ be a given profil in velocity. The linearized system around $(u_0,f_0)$ is
\beq \label{sys2}
\left\{
\begin{array}{lll}
\partial_t \widetilde u + A(u_0) \partial_x \widetilde u + \kappa  \partial_x \int_\rit \widetilde f(v) \phi(u_0,v) \, dv = 0, \\
\partial_t \widetilde f + v \partial_x \widetilde f + \Bigl( \nabla_u \psi(u_0) \cdot \partial_x \widetilde u \Bigr) \, \partial_v f_0 = 0.
\end{array}
\right.
\eeq
Then $\omega\in \mathbb C$ is an eigenvalue, with corresponding eigenfunction $u(t,x)=e^{i(kx-\omega t)}r$
where $0\neq k \in \mathbb R$ and $r\in \mathbb C^N$, if
\beq \label{dispp1}
\Bigl[ A(u_0) - \sigma \Bigr] r - \kappa \Bigl( \nabla_u \psi(u_0) \cdot r \Bigr)
\int_\rit {\phi(u_0,v) \partial_v f_0(v) \over v - \sigma} \, dv = 0, \qquad \sigma=\omega/k.
\eeq
Again this relation is established when $\Im \omega>0$ 
so that the interpretation of the integral is clear and it is compatible with the Fourier-Laplace transform. This formula is then extended by analyticity 
for all  $\omega$ in the  strip of analyticity of $f_0$.
When $\omega$ is real, we must change (\ref{dispp1}) into
\begin{equation} \label{dispp2}
\Bigl[ A(u_0) - \sigma \Bigr] r - \kappa \Bigl( \nabla_u \psi(u_0) \cdot u \Bigr) 
P.V.\int_\rit {\phi(u_0,v) \partial_v f_0(v) \over v - \sigma} \, dv
\end{equation}
$$
 - i \pi \kappa \Bigl( \nabla_u \psi(u_0) \cdot r \Bigr)  \phi(\omega) \partial_v f_0(\sigma) = 0.
$$
Let us now add the natural hypothesis that $A(u)$ comes from a strictly hyperbolic system, which means that $A(u)$ is the Jacobian matrix of some flux function and that $A(u) $ is diagonalizable in $\mathbb R$.
For symmetrizable systems,  the matrix  $A(u)$ can be chosen  symmetric
provided $u$ is a symmetrized variable
$$
A(u)=A(u)^t\in \mathcal M_{N}(\mathbb R).
$$
Then $A(u_0)$ has $N$ real eigenvalues $\sigma_j\in \mathbb R$ ($1 \le j \le N$) and $N$ real eigenvectors 
 $r_j\in \mathbb R^N$ ($1 \le j \le N$). 
The strict hyperbolicity hypothesis means that $\sigma_j\neq
\sigma_p$ for $j\neq p$. 

In (\ref{dispp2}) the eigenproblem depends on the parameter
$\kappa\in \mathbb R$, so it is natural to investigate the dependency
of the eigenvalues with respect to $\kappa$ in the vicinity of 
$\kappa=0$. It means that (\ref{dispp2}) is viewed
as a  eigenproblem with a pertubation of the form
\begin{equation} \label{pert}
\Bigl[ A(u_0) - \sigma \Bigr] r +
\kappa S(\sigma,r)=0
\end{equation}
where $S$ is non linear but analytic with respect to the eigenvalue $\sigma\in \mathbb C$ and is linear
with respect to the eigenvector $r\in \mathbb C^N$.
For small $\kappa$, it is natural as in \cite{kato} to admit that 
both $\omega_i(\kappa)$ and $r_j(\kappa)$ have a smooth
dependency with respect to $\kappa$.

Then we differentiate (\ref{pert}) with respect to $\kappa$, which leads to
$$
\Bigl[ A(u_0) - \sigma_j(\kappa) \Bigr] r_j'(\kappa)
- \sigma_j'(\kappa)  r_j(\kappa)+
\ S(\sigma_j(\kappa),r_j(\kappa))
+\kappa \left(\dots  \right)=0.
$$
For $\kappa=0$, we take the sesquilinear product against $r_j(0)$ and obtain
$$
\left(
r_j(0)\cdot \Bigl[ A(u_0) - \sigma_j(0) \Bigr] r_j'(0)
\right)
- \sigma_j'(0)  \|r_j(0) \|^2
+
\Bigl(
r_j(0)\cdot  S(\sigma_j(0),r_j(0))
\Bigr)=0.
$$
Due to the symmetry of the matrix, the first term vanishes 
 and one obtains
\begin{equation}
 \sigma_j'(0)=
 \frac{S(\sigma_j(0),r_j(0))}{\|r_j(0) \|^2}.
\end{equation}
In view of (\ref{dispp2}), this gives the variation of the imaginary part of the eigenvalue
$$
\left( \Im  \sigma_j\right) '(0)=-
\frac \pi {\|r_j(0) \|^2} \Bigl( \nabla_u \psi(u_0) \cdot r_j(0) \Bigr)  \left( \phi(u_0,\omega), r_j(0) \right) \partial_v f_0(\sigma_j(0)).
$$
We have thus proved the following proposition
which characterizes the linear stability of small kinetic perturbations of quasi-linear strictly hyperbolic systems systems of the form  $\partial_t u +A(u)\partial_xu=0$.

\begin{proposition} \label{pro:1}
A necessary condition for linear stability of solutions
of the system (\ref{sys1}) around $(u_0,f_0)$ 
in the vicinity of $\kappa=0$
is
$$
\Bigl( \nabla_u \psi(u_0) \cdot r_j \Bigr)  \left( \phi(u_0,\omega), r_j \right) \partial_v f_0(\sigma_j(0)) \geq 0
\quad 1\leq j \leq N
$$
where $(\omega_j,r_j)$ denotes any eigenpair of the symmetric matrix $A(u_0)$. 
\end{proposition}


Additionally if  $f_0$ is monotonic in $v$, if $\phi(u_0,\omega) \ne 0$ and if $\nabla_u \psi(u_0) \cdot v_j \ne 0$,
we see that, for small $\kappa$, $\omega_j(\kappa)$ can not be real: as in the previous example, 
the eigenvalues of the hyperbolic parts are ``expelled" from the real axis by their interaction with the kinetic part.
Note that $\nabla_u \psi(u_0) \cdot v_j = 0$ means that the eigenvector $v_j$ ``does not interact" with the kinetic part.


\subsubsection*{Acknowledgements}


The authors would like to thank C. Buet for many fruitful discussions.
The work of D. Bian is supported by NSF of China under the grant 12271032.


\end{document}